\documentclass{amsart}

\usepackage{graphicx}

\usepackage{amsthm}
\usepackage{amsmath}
\usepackage{amsfonts}

\usepackage{subcaption}

\usepackage{cite}

\author{Jean-Philippe Burelle}
\address{Jean-Philippe Burelle\\CNRS and Institut des Hautes \'Etudes Scientifiques\\35 route de Chartres\\91440 Bures-sur-Yvette\\France}
\email{jburelle@ihes.fr}

\author{Dominik Francoeur}
\address{Dominik Francoeur\\Universit\'e de Gen\`eve\\1205 Geneva\\Switzerland}
\email{dominik.francoeur@unige.ch}

\subjclass[2010]{20H10, and 53C50}

\thanks{The authors gratefully acknowledge support from the Natural Sciences and Engineering Research Council of Canada (NSERC), and U.S. National Science Foundation grants DMS 1107452, 1107263, 1107367 ``RNMS: Geometric Structures and Representation Varieties'' (the GEAR Network). The first author received funding for this project from the European Research Council (ERC) under the European Union’s Horizon 2020 research and innovation programme (ERC starting grant DiGGeS, grant agreement No 715982).}
\title{Foliations between crooked planes in $3$-dimensional Minkowski space}

\renewcommand{\vec}{\mathbf}
\newcommand{\ldot}[2]{#1 \cdot #2}
\newcommand{\lcross}[2]{\vec{#1} \times \vec{#2}}
\newcommand{\Min}{\mathrm{Min}}
\newcommand{\Ort}{\mathsf{O}}
\newcommand{\dif}{\mathrm{d}}

\newcommand{\bR}{\mathbb{R}}

\newcommand{\interior}{\mathrm{int}}

\newtheorem{thm}{Theorem}
\newtheorem{prop}{Proposition}

\theoremstyle{definition}
\newtheorem{defn}{Definition}[section]
\newtheorem{rmk}{Remark}[section]
\newtheorem{nota}{Notation}[section]

\begin{document}

\begin{abstract}
  We show that any two disjoint crooked planes in $\bR^3$ are leaves of a crooked foliation. This answers a question asked by Charette and Kim \cite{foliations}.
\end{abstract}

\maketitle

\section{Introduction}
In 1983, answering a question of Milnor\cite{Milnor}, Margulis constructed the first examples of nonabelian free groups which act freely and properly discontinuously on $\bR^3$ by affine transformations~\cite{Margulis}. In order to better understand these examples, Todd Drumm defined piecewise linear surfaces called \emph{crooked planes} which can bound fundamental domains for such actions~\cite{drummthesis}.

Crooked planes have proven to be very useful in the study of affine actions. Charette-Drumm-Goldman have used them in order to obtain a complete classification for free groups of rank two \cite{CDG1,CDG2,CDG3}. In particular, they show that every free and properly discontinuous affine action of a rank two free group on $\bR^3$ admits a fundamental domain bounded by finitely many crooked planes (the \emph{crooked plane conjecture}). A consequence of this is the \emph{tameness conjecture}, that the quotient of $\bR^3$ by one of these actions is homeomorphic to the interior of a compact manifold with boundary.

Building on this work, Danciger-Gu\'eritaud-Kassel showed in ~\cite{DGKArcComplex} that crooked planes have a natural interpretation in terms of the deformation theory of hyperbolic surfaces, and used this fact in order to prove the crooked plane conjecture in arbitrary rank, assuming that the linear part is convex cocompact in $\Ort(2,1)$.

One of the key aspects of the theory of crooked planes is their intersection properties. In particular, knowing when two crooked planes are disjoint is crucial. The \emph{Drumm-Goldman inequality} provides a necessary and sufficient criterion for two crooked planes to be disjoint ~\cite{DRUMM1999323}. This criterion was later expanded upon in ~\cite{halfspaces} and reinterpreted in terms of hyperbolic geometry in ~\cite{DGKArcComplex}.

As an application of the disjointness criterion, the first example of a \emph{crooked foliation}, a smooth $1$-parameter family of pairwise disjoint crooked planes, was given in ~\cite{halfspaces}. Charette-Kim \cite{foliations} investigated these foliations further and  gave necessary and sufficient criteria for a one-parameter family of crooked planes to foliate a subset of $\bR^3$. They ask the following question : given a pair of disjoint crooked planes in $\bR^3$, can the region between them be foliated by crooked planes? We answer this question in the affirmative.

\begin{thm}\label{thm:mainthm}
  Let $C,C'$ be a pair of disjoint crooked planes in $\bR^3$. Then, there is a \emph{crooked foliation}, that is, a smooth family of pairwise disjoint crooked planes $C_t$, $0\leq t\leq 1$ with  $C_0=C$ and $C_1=C'$.
\end{thm}

After recalling some definitions from the theory of crooked planes in Minkowski $3$-space in Section 2, we will prove the main theorem in Section 3.\\

We are thankful to the referee for insightful comments and for suggesting an elegant way to shorten the proof of the main theorem.

\section{Definitions}
\begin{defn}
Lorentzian $3$-space $\mathbb{R}^{2,1}$ is the real three dimensional vector space $\mathbb{R}^3$ endowed with the following symmetric bilinear form of signature $(2,1)$:
\[ \cdot : \mathbb{R}^3 \times \mathbb{R}^3 \rightarrow \mathbb{R}\]
\[ (\vec{u},\vec{v}) \mapsto u_1 v_1 + u_2 v_2 - u_3 v_3.\]
We fix the orientation given by the standard basis $e_1,e_2,e_3$ and we define the \emph{Lorentzian cross product}
\[\vec{u}\times \vec{v} = (u_2 v_3 - u_3 v_2, u_3 v_1 - u_1 v_3, u_2 v_1 - u_1 v_2)\in\bR^{2,1},\]
for $\vec{u},\vec{v}\in\bR^{2,1}$.
\end{defn}

A \emph{null frame} of $\bR^{2,1}$ is a positively oriented basis $\vec{u},\vec{u}',\vec{u}''$ such that $\ldot{\vec{u}}{\vec{u}}=1$, $\ldot{\vec{u}'}{\vec{u}''}=-1$ and all other products between the three vectors vanish.

\begin{nota}
  Any unit spacelike vector $\vec{u}$ can be extended to a null frame. This frame is unique up to scaling $\vec{u}'$ and $\vec{u}''$ by inverse scalars. As normalization we will choose $\vec{u}'$ and $\vec{u}''$ so that their third coordinates are positive and equal. Given $\vec{u}$, we will denote these two null vectors by $\vec{u}^-$ and $\vec{u}^+$, respectively.
\end{nota}

We will denote by $\Min$ the pseudo-Euclidean affine space which is modeled on the vector space $\mathbb{R}^{2,1}$. In other words, $\Min$ is a topological space on which $\mathbb{R}^{2,1}$ acts simply transitively by homeomorphisms. For $\vec{v}\in \mathbb{R}^{2,1}$ and $p\in\Min$, we denote this action by $\vec{v}(p) = p + \vec{v}$. If $q=p+\vec{v}$, we will also write $q-p=\vec{v}$. A choice of origin $o\in \Min$ identifies $\Min$ with $\mathbb{R}^{2,1}$ via the map $\vec{v}\mapsto o + \vec{v}$.

\begin{figure}[h]
\begin{minipage}[b]{.7\linewidth}
  \centering\includegraphics[width=.7\linewidth]{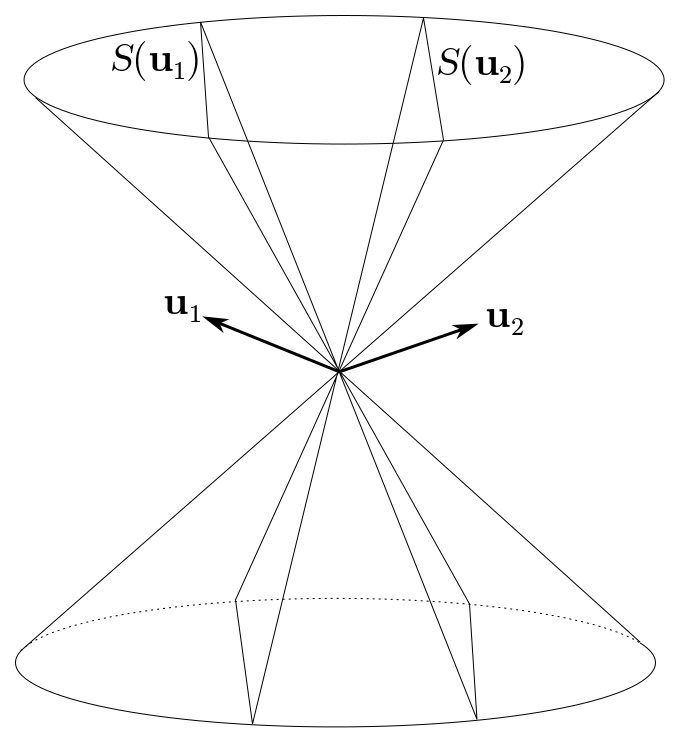}
  \subcaption{A pair of consistently oriented unit spacelike vectors $\vec{u}_1,\vec{u}_2$ and the corresponding stems.}
\end{minipage}
\begin{minipage}[b]{.7\linewidth}
  \centering\includegraphics[width=.9\linewidth]{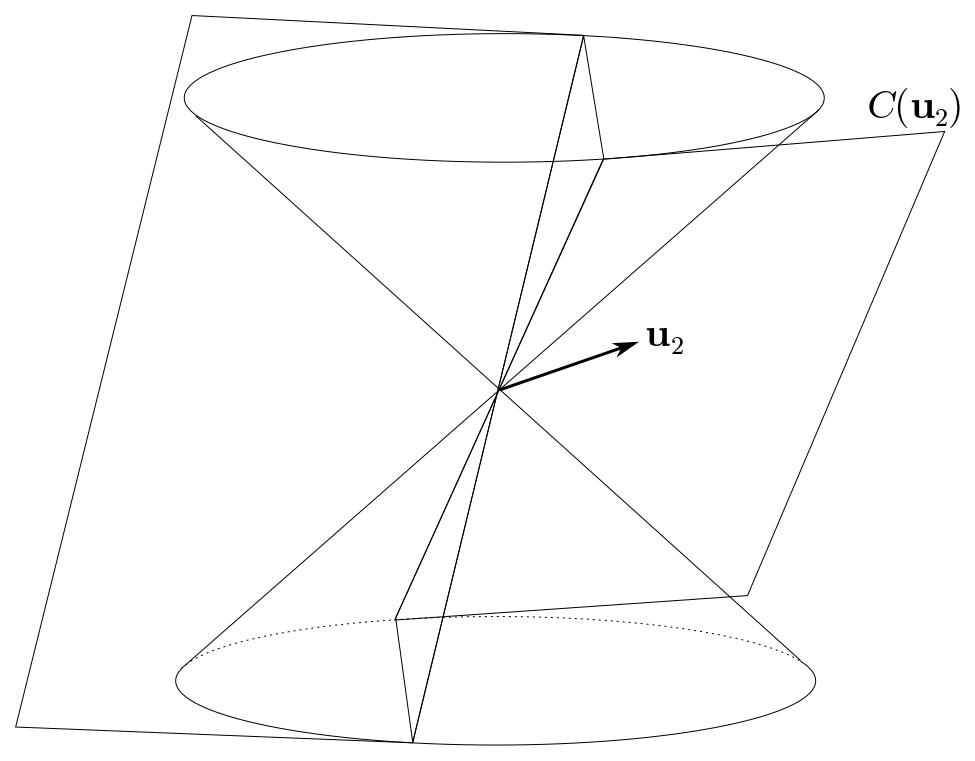}
\subcaption{The linear crooked plane $C(\vec{u}_2)$.}
\end{minipage}
\caption{Consistent orientations, stems and crooked planes.}
\label{fig:consistorient}
\end{figure}

We now recall the definition of a crooked plane. First, we define a \emph{stem}, which will be one of the three linear pieces of a crooked plane.
\begin{defn}
Let $\vec{u}\in \mathbb{R}^{2,1}$ be a unit spacelike vector. The \emph{stem} $S(\vec{u})$ is the set of causal vectors orthogonal to $\vec{u}$ :
\[S(\vec{u}) = \{ \vec{v} \in \mathbb{R}^{2,1} ~|~ \ldot{\vec{u}}{\vec{v}}=0 \text{ and } \ldot{\vec{v}}{\vec{v}} \leq 0 \}.\]
A stem is the union of two opposite closed quadrants (see Figure \ref{fig:consistorient}).
\end{defn}
\begin{defn}
Let $\vec{u}\in \mathbb{R}^{2,1}$ be a unit spacelike vector. The \emph{linear crooked plane} $C(\vec{u})$ is the piecewise linear surface defined by:
\[ C(\vec{u}) := \{ \vec{v} \in \mathbb{R}^{2,1} ~|~ \lcross{v}{w} = k \vec{w} \textrm{ for some }  \vec{w}\in S(\vec{u}) \text{ and } k\in \mathbb{R}_{\geq 0}\}.\]
\end{defn}
From this definition, we see that $S(\vec{u})\subset C(\vec{u})$ since $\vec{v}\times\vec{v}=0$ for all $\vec{v}\in \bR^{2,1}$. The complement of the stem $C(\vec{u})-S(\vec{u})$ has two connected components which are called the \emph{wings} of the crooked plane. Each wing is a half-plane on which the Lorentzian bilinear form is degenerate, attached to the stem along its boundary (See Fig. \ref{fig:consistorient}). Note that $C(\vec{u})=C(-\vec{u})$.

\begin{defn}
Let $p\in \Min$ and $\vec{u}\in \mathbb{R}^{2,1}$ unit spacelike. The \emph{crooked plane} $C(p,\vec{u})$ is the set $p + C(\vec{u}) \subset \Min$. The vector $\vec{u}$ is called a \emph{directing vector} of the crooked plane, and $p$ its \emph{vertex}.
\end{defn}

In order to formally state the disjointness criteria from \cite{DRUMM1999323,foliations}, we need a normalization for pairs of unit spacelike vectors.
\begin{defn}
  Two unit spacelike vectors $\vec{u}_1,\vec{u}_2 \in \mathbb{R}^{2,1}$ are \emph{consistently oriented} if
  \begin{itemize}
    \item $\ldot{\vec{u}_1}{\vec{u}_2}\leq -1$, and
    \item $\ldot{\vec{u}_i}{\vec{u}_j}^\pm\leq 0$ for $1\leq i,j\leq 2$.
  \end{itemize}
\end{defn}

Two consistently oriented unit spacelike vectors $\vec{u},\vec{u}'$ are called \emph{ultraparallel} if $\ldot{\vec{u}}{\vec{u}'}<-1$. They are called \emph{asymptotic} if $\ldot{\vec{u}}{\vec{u}'}=-1$ and $\vec{u}'\neq -\vec{u}$. Intersecting $\vec{u}^\perp$ and $\vec{u'}^\perp$ with the hyperboloid model of the hyperbolic plane defines a pair of hyperbolic geodesics, and the terminology comes from the relative position of these geodesics. Choosing one of the unit vectors $\pm\vec{u}$ endows the geodesic in the hyperboloid model of $\mathbb{H}^2$ defined by $\vec{u}^\perp$ with a transverse orientation. Two unit spacelike vectors are consistently oriented when the corresponding transversely oriented geodesics are disjoint with transverse orientations pointing away from each other (see Figure \ref{fig:consistorient}).

Whenever there exists a choice of directing vectors $\vec{u},\vec{u}'$ which are consistently oriented, we will also call a pair of crooked planes $C(p,\vec{u}),C(p',\vec{u}')$ ultraparallel or asymptotic accordingly.

We will use two disjointness criteria for crooked planes, one for pairs of crooked planes and one for foliations. Both depend on the following notion:
\begin{defn}
  The \emph{stem quadrant} associated to a unit spacelike vector $\vec{u}$ is the set
  \[\mathsf{V}(\vec{u}) := \{a \vec{u}^- - b\vec{u}^+ : a,b\geq 0\}\backslash\{0\}.\]
  Note that $\mathsf{V}(-\vec{u})=-\mathsf{V}(\vec{u})$.
\end{defn}

The following disjointness criterion is a restatement of the \emph{Drumm-Goldman inequality} ~\cite{DRUMM1999323}.
\begin{thm}[Burelle-Charette-Goldman \cite{halfspaces}]\label{DrummgoldmanDisjoint}
  Let $C=C(p,\vec{u})$, $C'=C(p',\vec{u}')$ be crooked planes and assume that $\vec{u},\vec{u}'$ are consistently oriented. Then, $C$ and $C'$ are disjoint if and only if
  \[p'-p \in \mathsf{A(\vec{u},\vec{u}')}:=\interior(\mathsf{V}(\vec{u}')-\mathsf{V}(\vec{u})).\]
\end{thm}

\begin{rmk}
  It is also shown in \cite{DRUMM1999323} that if there is no choice of sign for $\vec{u},\vec{u}'$ making them consistently oriented, then $C(p,\vec{u})$ and $C(p',\vec{u}')$ necessarily intersect. Therefore, the above theorem is a characterization of disjoint crooked planes.
\end{rmk}

We will use the following straightforward consequence of the \emph{Charette-Kim criterion} for crooked foliations (foliations of $\bR^{2,1}$ by crooked planes) :
\begin{thm}[Charette-Kim \cite{foliations}]\label{charettekimUltrap}
  Let $(\vec{u}_t)_{t\in \mathbb{R}}$ be a path of pairwise ultraparallel or asymptotic unit spacelike vectors such that $-\vec{u}_t,\vec{u}_s$ are consistently oriented for all $t<s$. Suppose $(p_t)_{t\in \mathbb{R}}$, is a regular curve such that for every $t\in\mathbb{R}$,
  \[\dot{p_t} \in \interior(\mathsf{V}(\vec{u}_t)).\]
  Then, $C(p_t,\vec{u}_t)$ is a crooked foliation.
\end{thm}


\section{Foliations between crooked planes}

We now prove Theorem \ref{thm:mainthm} : there exists a crooked foliation containing any pair of disjoint crooked planes. The theorem is a consequence of the following stronger result :

\begin{prop}
  Let $(\vec{u}_t)_{t\in[0,1]}$ be a smooth path of unit spacelike vectors which are pairwise ultraparallel or asymptotic. Let $p_0,p_1\in \Min$ such that $C(p_0,\vec{u}_0)$ and $C(p_1,\vec{u}_1)$ are disjoint crooked planes. Then, there exists a path $(p_t)_{t\in [0,1]}$ starting at $p_0$ and ending at $p_1$ such that $C(p_t,\vec{u}_t)$ is a smooth crooked foliation.

  \begin{proof}

Since we assume that $\vec{u}_s$ are pairwise ultraparallel or asymptotic, we have that $\ldot{\vec{u}_t}{\vec{u}_s}\geq 1$ for all $t\leq s$. Changing the path $\vec{u}_s$ to $-\vec{u}_s$ if needed (both paths define the same linear crooked planes) we may also assume that $-\vec{u}_t,\vec{u}_s$ are consistently oriented for all $t<s$.

For any pair of smooth functions $f,g : [0,1]\rightarrow \bR^{>0}$, define
\[\vec{v}_{f,g}(s) := f(s)\vec{u}_s^- - g(s)\vec{u}_s^+.\]
Then, the path of vertices $p_{f,g}(t) := p_0 + \int_0^t \vec{v}_{f,g}(s)\,\dif s$ satisfies the hypotheses of Theorem \ref{charettekimUltrap} since its derivative
\[\dot{p}_{f,g}(t) = \vec{v}_{f,g}(t)\]
lies in the interior of $\mathsf{V}(\vec{u}_t)$ by definition.

Let $\mathsf{D}$ denote collection of displacement vectors $p_{f,g}(1)-p_0$ :
\[\mathsf{D} = \left\{\left. \int_0^1 \vec{v}_{f,g}(s)\,\dif s ~\right|~ f,g :[0,1]\rightarrow \bR^{>0}\right\}.\]
Then $\mathsf{D}$ is a convex cone since $k\vec{v}_{f,g}=\vec{v}_{kf,kg}$ for $k\in \bR^{>0}$ and $\vec{v}_{f_1,g_1} + \vec{v}_{f_2,g_2} = \vec{v}_{f_1+f_2,g_1+g_2}$. Moreover, since by Theorem \ref{charettekimUltrap} the crooked planes $C(p_{f,g}(t),\vec{u}_t)$ define crooked foliations, the initial and final crooked planes are disjoint and so $\mathsf{D}\subset\mathsf{A}(-\vec{u}_0,\vec{u}_1)$. Since the cone $\mathsf{A}(-\vec{u}_0,\vec{u}_1)$ is the interior of the convex hull of the four rays generated by $\vec{u}_0^-,-\vec{u}_0^+,\vec{u}_1^-,-\vec{u}_1^+$, to show equality of the cones it suffices to show that these rays can be approximated by vectors in $\mathsf{D}$.

Consider the sequences $f_n(s) = ne^{-n s}$ and $g_n(s)=e^{-n}$. Integrating by parts we get
\[\int_0^1 f_n(s)\vec{u}_s^-\,\dif s = \vec{u}_0^- - e^{-n}\vec{u}_1^- + \int_0^1 e^{-n s}\dot{\vec{u}}_s^- \,\dif s.\]
Therefore, as $\vec{u}_s$ is smooth and so $\vec{u}^+_s$ and $ \dot{\vec{u}}^-_s$ are bounded on $[0,1]$,
\[\lim_{n\rightarrow \infty} \int_0^1\vec{v}_{f_n,g_n}(s)\,\dif s = \vec{u}^-_0.\]
We conclude that $\mathsf{D}$ contains vectors arbitrarily close to the ray $\bR^{>0}\vec{u}_0^-$.

Similarly, if $f$ is concentrated near $s=1$ and $g$ is small we can approximate the ray $\vec{u}_1^-$, and exchanging the roles of $f$ and $g$ we approximate the other two rays on the boundary of the convex cone $\mathsf{A}(\vec{u}_0,\vec{u}_1)$.
\end{proof}
\end{prop}

The previous proposition has the following interpretation : given any geodesic foliation $\mathcal{F}$ of the region between two geodesics $\ell_0,\ell_1$ of $\mathbb{H}^2$ and basepoints $p_0,p_1\in \Min$ such that the crooked planes with vertices $p_i$ and stems corresponding to $\ell_i$ are disjoint, $\mathcal{F}$ can be lifted to a foliation by crooked planes of the region between the crooked planes.

\bibliography{mybib}{}
\bibliographystyle{plain}
\end{document}